\documentstyle[12pt]{article} \textwidth 6.1 in \textheight 9.5 in

\title{\bf Infinite multiplicity for inhomogeneous supercritical problem in entire space
 \thanks { This work was supported by the the Natural Science Foundation of
China (No:10901047)  }}
\author{Baishun Lai,  Zhihao  Ge\\
\small {\it School of Mathematics and information sciences, Henan  University}\\
\small {\it Kaifeng,  Henan 475001, China}\\
} \baselineskip 0.2in
\date{}
\begin{document}
\maketitle
\begin{center}
\begin{minipage}{130mm}
{\small {\bf Abstract}\vskip0.1in
 In this paper, we will prove the existence of infinitely many  positive solutions to the following
supercritical problem by using the Liapunov-Schmidt reduction method
and asymptotic analysis:
\begin{eqnarray*}
 \left\{
\begin{array}{ll}\Delta u + u^{p}+f(x)=0,\ \ u>0\ \ \mbox{in}\
R^{n},\\
 \lim_{|x|\to\infty}u(x)\to 0.
\end{array}
\right.
\end{eqnarray*}

 }
\end{minipage}
\end{center}
\baselineskip 0.2in \vskip 0.2in {\bf 1. Introduction and statement of the main results} %
\setcounter{section}{1} \setcounter{equation}{0}\vskip 0.2in
   The purpose of this paper is to establish the existence of
   infinitely many positive solutions to the following inhomogeneous
   problem
   \begin{equation}\left\{\begin{array}{ll}
   \Delta u + u^{p}+f(x)=0,\\
   u>0\ \ \mbox{in}\ R^{n},\ \lim_{|x|\to\infty}u(x)=0.
   \end{array}\right.
\end{equation}
Where $n\geq 3,p>\frac{n+2}{n-2},
\Delta=\sum_{i=1}^{n}(\partial^{2}/\partial x_{i}^{2})$ is the
Laplacian operator and $f(x)\in C^{0,\alpha}_{loc}(R^{n})$ with
$f\geq 0$ everywhere in $R^{n}$, $f\not \equiv0$.\vskip0.1in

Inhomogeneous second-order elliptic equations defined in entire
space arise naturally in probability theory in the study of
stochastic processes. the eqns. (1.1) in particular appeared recently in a
paper by Tzong-Yow Lee \cite{Fowler} establishing limit theorems for
super-brownian motion. In that paper, existence results for the eqns. (1.1)
were obtained in the case where the inhomogeneous terms is compactly
supported and is dominated by a function of the form
$\frac{C}{(1+|x|)^{(n-2)p}} $ where $C(n,p)>0$ is sufficiently
small. In addition to Lee [9], (1.1) has been studied by Pokhozhaev
[11] and Egnell and Kay [10]. Pokhozhaev obtained radial solutions
when the inhomogeneous term is radially symmetric about the origin
and satisfies certain integrability conditions. Also, Egenll and Kay
worked on an equation similar to that of (1.1) but have a small
positive parameter $\varepsilon$ as a coefficient to inhomogeneous
term $f(x)$. In a recent paper [1], existence results were obtained
for a large class of functions for the inhomogeneous term, i.e.,
$$0\leq f(x)\leq \frac{p-1}{[p(1+|x|^{2}]^{\frac{p}{p-1}}}L^{p}$$
for all $x\in R^{n}$, where
$$L=\Big[\frac{2}{p-1}(n-2-\frac{2}{p-1})\Big]^{\frac{1}{p-1}}.$$

Using sub-super solution method, Bae and Ni established the
following infinite multiplicity result for the equation
\begin{eqnarray*}
 \left\{
\begin{array}{ll}\Delta u + u^{p}+\mu f(x)=0,\ \ u>0\ \ \mbox{in}\
R^{n},\\
 \lim_{|x|\to\infty}u(x)\to 0.
\end{array}
\right.
\end{eqnarray*} where
$\mu>0$ is a parameter.\vskip 0.1in

{\bf Theorem A$^{[1]}$:} (i) Let $p>p_{c}$. Suppose that near
$\infty$
$$ \max(\pm f(x), 0)\leq |x|^{-q_{\pm}},$$ where
$q_{+}>n-\lambda_{2}$ and $q_{-}>n-\lambda_{2}-\frac{2}{p-1}$. Then,
there exists $\mu_{*}>0$ such that for every $\mu \in (0,\mu_{*})$,
equation (1.2) possesses infinitely many solutions with asymptotic
behavior $L|x|^{-m}$ at $\infty$

(ii) Let $p=p_{c}$. Then, the conclusion in (i) holds if we assume
in addition that either $f$ has a compact support in $R^{n}$ or $f$
does not change sign in $R^{n}$, where
\[
p_{c}=\left\{
\begin{array}{ll}
\frac{(n-2)^{2}-4n+4\sqrt{n^{2}-(n-2)^{2}}}{(n-2)(n-10)}\ \ \ \ \ \
\ n>10,
&  \\
\infty \ \ \ \ \ \ \ \ \ \ \ \ \ \ \ \ \ \ \ \ \ \ \ \ \ \ \ \ \ \ \
\ \ \ 3\leq n\leq 10. &
\end{array}%
\right.
\]
and
$$\lambda_{2}=\frac{(n-2-\frac{4}{p-1})+\sqrt{(n-2-\frac{4}{p-1})^{2}-8(n-2-\frac{2}{p-1})}}{2}$$

Besides, for a more general case (i.e., $u^{p}+\mu f$ is replaced by
$K(x)u^{p}+\mu f$), similar results are obtained, see [2, 8, 12].
However, these results are established require $p>p_{c}, n\geq 11$.
But the case of $n\leq 10$ and $\frac{n+2}{n-2}<p<p_{c}$ is {\bf
open}.
 Note that in this case, the
method of sub-super solutions breaks down. In this paper, under
reasonable conditions on $f$, we establish that when
$p>\frac{n+1}{n-3}$, (1.1) has a continuum of solutions.\vskip 0.1in

The main difficulties in establishing existence of (1.1), in
addition to the noncompactness of the domain and the presence of an
inhomogeneous term, are the lack of (local) sobolev embedding
suitably fit to a weak formulation of this problem. So direct tools
of the calculus of variation are not appropriate for (1.1).\vskip
0.1in

Instead of using  sub-super solution method (which limits the
applicability on the exponent $p$), we use asymptotic analysis and
Liapunov-Schmidt reduction method to prove Theorem 1.  This is based
on the construction of a sufficiently good approximation solution.
It is well known that the problem
\begin{equation} \Delta w +w^{p}=0 \ \ \mbox{in}\ R^{n}\end{equation}
where $p\geq\frac{n+2}{n-2}$, possesses a positive radially
symmetric solution $w(|x|)$ which reduces to the equation
\begin{equation}w''+\frac{n-1}{r}w'+w^{p}=0.\end{equation} This
equation can be analyzed through phase plane analysis after a
transformation introduced  by Fowler [9] in 1931:
$v(s)=r^{\frac{2}{p-1}}w(r), r=e^{s},$ which transforms equation
(1.3) into the autonomous ODE \begin{equation}v''+\alpha v'-\beta v
+v^{p}=0\end{equation} where $$\alpha=n-2-\frac{4}{p-1},
\beta=\frac{2}{p-1}(n-2-\frac{2}{p-1}).$$ Since $\alpha,\beta$ are
positive for $p>\frac{n+2}{n-2}$, the Hamiltonian energy
$$E(v)=\frac{1}{2}v'^{2}+\frac{1}{p+1}v^{p+1}-\frac{\beta}{2}v^{2}$$
strictly decrease along trajectories. Using this it is easy to see
the existence of a heteroclinic orbit which connects the equilibria
(0,0) and $(0,\beta^{\frac{1}{p-1}})$ in the phase plane $(v,v')$,
which corresponding respectively to a saddle point and an attractor.
A solution of (1.4) corresponding this orbit satisfies
$v(-\infty)=0, v(+\infty)=\beta^{\frac{1}{p-1}}$ and
$w(r)=r^{-\frac{2}{p-1}}v(\log r)$ solves (1.3) and is bounded at
$r=0$. Then all radial solutions of (1.3) defined in all $R^{n}$
have the form
$$w_{\lambda}(x):=\lambda^{\frac{2}{p-1}}w(\lambda|x|), \ \
\lambda>0.$$

We denote in what follows by $w(x)$ the unique positive radial
solution $$ \Delta w+w^{p}=0\ \ \mbox{in}\ \ R^{n}, \ \ \ w(0)=1.$$

 At main
order one has $$w(r)\sim L r^{-\frac{2}{p-1}}\ \ \mbox{as} \ r\to
\infty,$$ which implies that this behavior is actually common to all
solutions $w_{\lambda}(x)$. \vskip 0.1in

By the change of variable
$\lambda^{-\frac{2}{p-1}}u(\frac{x}{\lambda})$ the (1.1) becomes
\begin{equation}
 \Delta u +u^{p}+f_{\lambda}(x)=0, \ \ \mbox{in}\ \ R^{n} \ \  u>0,
\ \ \lim_{|x|\to+\infty}u(x)=0, \label{eq1.5}
\end{equation}
where
$f_{\lambda}(x)=\lambda^{-\frac{2p}{p-1}}f(\frac{x}{\lambda})$. In fact,
if $f$ is assumed to satisfy the asymptotic behavior
$$f(x)=o(|x|^{-\frac{2p}{p-1}})\ \  \mbox{as}\ \ |x|\to+\infty,$$
then we observe that away from the origin $f_{\lambda}(x) \to 0$ as
$\lambda \to 0$. Thus (1.1) may be regarded, away from the origin,
as small perturbations of problem (1.2) when $\lambda>0$ is
sufficiently small. From (\ref{eq1.5}), we find that (1.1) "hide" a
parameter which indexes a continuum of solutions which
asymptotically vanish over compact sets.\vskip0.1in

Our main result is as follows:\vskip0.1in

 {\bf Theorem 1.} Let $p>\frac{n+1}{n-3}, n\geq 4$. Assume that
$$f(x)=\eta(x)f_{1}(x),$$ where $\eta(x)\in C_{0}^{\infty}(R^{n}),
0\leq \eta(x)\leq 1,$ $$ \eta(x)=0\ \ \mbox{for}\ \ |x|\leq R_{1},
\eta(x)=1 \ \ \mbox{for}\ |x|\geq R_{1}+1,
$$$(R_{1}
> 0$ fixed large enough) and $f_{1}(x)\in C^{0,\alpha}(R^{n}), f_{1}(x)<|x|^{-\mu}$ with
$\mu>2+\frac{2}{p-1}$. Then problem (1.1) has a continuum of
solutions $u_{\lambda}(x)$ (parameterized by $\lambda\leq
\lambda_{0}$, where $\lambda_{0}$ is a fixed number ) such that
$$\lim_{\lambda \to 0} u_{\lambda}(x)=0$$ uniformly in
$R^{n}\setminus \{0\}$. The same results holds when
$\frac{n+2}{n-2}<p\leq \frac{n+1}{n-3}$ provided that $f$ is
symmetric with respect to $n$ coordinate axis, namely
$$f(x_{1},...,x_{i},...,x_{n})=f(x_{1},...-x_{i},...x_{n}),\ \ \mbox{for\ all}\ i=1,...,n.$$

The idea is to consider $w(x)$ as an approximation for a solution of
(1.1), provided that $\lambda>0$ is chosen small enough. To this
end, we need to study the solvability of the operator
$\Delta+pw^{p-1}$ in suitable weighted Sobolev space. Recently, this
issue has been studied in Davila-del Pino-Musso \cite{4Pino-Wei} and
Davila-Pino-Musso-Wei [5].

Throughout the paper, the symbol $C$ denotes always a positive
constant independent of $\lambda$, which could be changed from one
line another. Denote $A\sim B$ if and only if there exist two
positive numbers $a,b$ such that $aA\leq B\leq bA$.\vskip0.2in

{\bf 2. The solvability of linearized operator $\Delta +pw^{p-1}$}
\setcounter{section}{2} \setcounter{equation}{0} \vskip0.1in

Our main concern in this section is to state the results concerning
the existence of solution in certain weighted spaces for
\begin{equation} \Delta \phi+pw^{p-1}\phi=h\ \ \mbox{in}\ R^{n},\end{equation} where
$w$ is the radial solution to (1.2) and $h$ is a known function
having a specific decay at infinity. We are looking for a solution
to (2.1) that is turn out to be a perturbation of $w$, it is rather
natural to require that it has a decay at most the same as that of
$w$, namely $O(|x|^{-\frac{2}{p-1}})$ as $|x|\to \infty$. Of course
we would also like $\phi$ be bounded on compact sets. As a result,
we shall assume that $h$ behaves like this but with two powers
subtracted, that is, $h=O(|x|^{-\frac{2}{p-1}-2}$) at infinity.

Now, define some weighted $L^{\infty}$ norms as follows (adopted
from [4,5]): $$\|\phi\|_{*}=\sup_{|x|\leq
1}|x|^{\sigma}|\phi(x)|+\sup_{|x|\geq
1}|x|^{\frac{2}{p-1}}|\phi(x)|,$$ and $$\|h(x)\|_{**}=\sup_{|x|\leq
1}|x|^{2+\sigma}|h(x)|+\sup_{|x|\geq
1}|x|^{2+\frac{2}{p-1}}|h(x)|,$$ where $\sigma>0$ will be fixed
later as needed.

The following lemmas and remarks on the solvability are due to
Davila-del Pina-Musso [4] and Davila-del Pino-Musso-Wei [5]:
\vskip0.2in

{\bf Lemma 2.1.} Assume that $p>\frac{n+1}{n-3},n\geq 4.$ For
$0<\sigma<n-2$ there exists a constant $C>0$ such that for any $h$
with $\|h(x)\|_{**}<\infty$, equation (2.1) has a solution
$\phi=T(h)$ such that $T$ defines a linear map and
$$\|T(h)\|_{*}\leq C \|h\|_{**}.$$ \vskip 0.2in

%Although  Lemma 2.1 is proven in [4, 5],
For the sake of completeness, we give the main idea of the proof of  Lemma 2.1  as follows (for the details, see [4, 5]):
%we summarize the main points of the argument

Let $\Theta_{k}, k\geq 0$ be the eigenfunction of the
Laplace-Beltrami operator $-\Delta_{s^{n-1}}$ on the sphere
$S^{n-1}$ with eigenvalues $\lambda_{k}$ repeated according to their
multiplicity, normalized so that they constitute an orthonormal
system in $L^{2}(S^{n-1})$. We let $\Theta_{0}$ be a positive
constant, associated to the eigenvalue 0 and $\Theta_{i}, 1\leq
i\leq n$ is an appropriate multiple of $\frac{x_{i}}{|x|}$ which has
eigenvalue $\lambda_{i}=n-1, 1\leq i\leq n$. we repeat eigenvalues
according to their multiplicity and we arrange them in an
non-decreasing sequence. We recall that the set of eigenvalues is
given by $\{j(n-2+j)| j\geq 0\}$. We write $h$ as
\begin{equation} h(x)=\Sigma_{k=0}^{k=+\infty}h_{k}\Theta_{k}(\theta)\end{equation}
and look for a solution $\phi$ to (2.1) in the form
\begin{equation} \phi(x)=\Sigma_{k=0}^{k=+\infty}\phi_{k}\Theta_{k}(\theta).\end{equation}
Then
\begin{equation}
\phi_{k}''+\frac{n-1}{r}\phi'+(pw^{p-1}-\frac{\lambda_{k}}{r^{2}})\phi=h_{k}
\end{equation}
Equation (2.4) can be solved for each $k$ separately:

{\bf a:} If $k=0$ and $p>\frac{n+2}{n-2}, \|h_{0}\|_{**}<+\infty$
then (2.4) has a solution $\phi_{0}$ which depends linearly on
$h_{0}$ and satisfies
\begin{equation}
\|\phi_{0}\|_{*}\leq C \|h_{0}\|_{**}.
\end{equation}Indeed in this case this solution is defined using the
variation of parameters formula
\begin{equation}\phi_{0}(r):=z_{1,0}(r)\int_{1}^{r}z_{2,0}h_{0}s^{n-1}ds-
z_{2,0}(r)\int_{0}^{r}z_{1,0}h_{0}s^{n-1}ds,\end{equation} where
$z_{1,0},z_{2,0}$ are two special linearly independent solution to
(2.4) with $k=0$ and $h_{0}=0$. More precisely, we take
$z_{1,0}=rw'+\frac{2}{p-1}w$ and $z_{2,0}$ a linearly independent
solution. Linearization shows that
$$z_{j,0}(r)=O(r^{-\frac{n-2}{2}})\ \ \mbox{as}\ \ r\to +\infty,\ \
j=1,2,$$ while $$z_{2,0}\sim r^{2-n}\ \ \mbox{near}\ \ r=0.$$ Using
this and (2.6), we can easily get estimate (2.5)

{\bf b:}  If $k=1, n\geq 4$ and $p>\frac{n+1}{n-3},
\|h_{1}\|_{**}<+\infty$, then we have
\begin{equation}
\|\phi_{1}\|_{*}\leq C \|h_{1}\|_{**}.
\end{equation} In this case, we have that the positive function
$z_{1}:=-w'(r)>0 $ in $(0,+\infty)$ solves (2.4) with $k=1$ and
$h_{1}=0$. Using this, we then define $\phi_{1}(r)$ as
\begin{equation}
\phi_{1}(r)=-z_{1}(r)\int_{1}^{r}z_{1}(s)^{-2}s^{1-n}ds\int_{0}^{s}z_{1}(\tau)h_{1}(\tau)\tau^{n-1}.
\end{equation}
Using this formula and by a simple computer, estimate (2.7) is
easily obtained.

{\bf c:} Let $k\geq 2$ and $p>\frac{n+2}{n-2}$. If $
\|h_{k}\|_{**}<+\infty$ (2.4) has a unique solution $\phi_{k}$ with
$\|\phi_{k}\|_{*}<+\infty$ and there exists $C_{k}>0$ such that
\begin{equation}
\|\phi_{k}\|_{*}\leq C \|h_{k}\|_{**}
\end{equation}
this case is simpler because the operator
$$L_{k}\phi=\phi''+\frac{n-1}{r}\phi'+(pw^{p-1}-\frac{\lambda_{k}}{r^{2}})\phi$$
satisfies the maximum principle in any interval of the form
$(\delta,\frac{1}{\delta}), \delta>0$. Indeed let $z=-w'$, so that
$z>0$ in $(0,+\infty)$ and it is a supersolution, because
$$L_{k}z=\frac{n-1-\lambda_{k}}{r^{2}}z<0\ \ \mbox{in}\ \
(0,+\infty),$$ since $\lambda_{k}\geq 2n$ for $k\geq n$. We
construct a supersolution $\psi$ of the form
$$\psi=C_{1}z+v,\ \ v(r)=\frac{1}{r^{\sigma}+r^{\frac{2}{p-1}}}.$$
Choosing $C_{1}$ sufficiently large, we can check that
$$L_{k}\psi\leq-c\min(r^{-\sigma-2},r^{-\frac{2}{p-1}-2})\ \
\mbox{in}\  (0,+\infty).$$ Using this, we can easily obtain (2.9).

The previous construction and (2.5), (2.7) and (2.9) imply that
given  an integer $m>0$, if $\|h\|<+\infty$ satisfies $h_{k}\equiv
0\ \ \forall k\geq m$ then there exists a solution $\phi$ to (3.1)
that depends linearly with respect to $h$ and moreover
$$\|\phi\|_{*}\leq C_{m} \|h\|_{**}$$ where $C_{m}$ may depend on
only $m$. Then by using a blow up argument, we can show that $C_{m}$
can be chosen independently of $m$. For a detailed proof, we refer
the interested readers to [5,6]. \ \  $\Box$\vskip0.1in

{\bf Remark 2.1.} From the above proof, we know the operator $T$ in
Lemma 2.1 are constructed \textquotedblleft by hand" by decomposing
$h$ and $\phi$ into suns of spherical harmonics where the
coefficients are radial functions. The nice property is of course
that $w$ is radial, the problem decouples into an infinite
collection of ODEs.\vskip0.1in

{\bf Remark 2.2.}  If $p\geq\frac{n+2}{n-2}$, linearized operator
$\Delta+pw^{p-1}$ has a kernel, i.e., $\mbox{span}\{\frac{\partial
w}{\partial x_{1}}, i=1,2,....n\}$, in general Sobolev space.
However, under suitable weighted Sobolev space, the linearized
operator $\Delta+pw^{p-1}$ is invertible, i.e., the kernel is
0.\vskip0.1in

{\bf 3. The proof of Theorem 1} \setcounter{section}{3}
\setcounter{equation}{0} \vskip0.1in

Let $p>\frac{n+1}{n-3}$, we will prove Theorem 1 in this section.
The main idea is to use Proposition 2.1 and a contraction mapping
principle.

We look for a solution of (1.5) of the form $u=w+\phi$, which yields
the following equation for $\phi$ $$\Delta
\phi+pw^{p-1}\phi=N(\phi)-f_{\lambda}(x),$$ where
\begin{equation} N(\phi)=-(w+\phi)^{p}+w^{p}+pw^{p-1}\phi. \end{equation}

Using the operator $T$ defined in Proposition 2.1, we are led to
solving the fixed point problem
\begin{equation}\phi=T(N(\phi)-f_{\lambda}(x)).\end{equation}

Firstly let us estimate $\|N(\phi)\|_{**,\lambda}$ depending on
whether $p\geq2$ or $p<2$.

Case $p\geq 2$. In this case, we observe that
$$|N(\phi)|\leq C(w^{p-2}\phi^{2}+|\phi|^{p}).$$ Let us
work with $0<\sigma\leq\frac{2}{p-1}$. Since
$$|\phi(x)|\leq
C|x|^{-\frac{2}{p-1}}\|\phi\|_{*},\ \ \mbox{for\ all}\ |x|\geq 1,$$
and
$$w(x)\leq C(1+|x|)^{-\frac{2}{p-1}},\ \mbox{for\ all}\
x\in R^{n},$$ so we have on one hand \begin{equation} \sup_{|x|\geq
1}|x|^{2+\frac{2}{p-1}}w^{p-2}|\phi|^{2} \leq
C\|\phi\|^{2}.\end{equation} On the other hand,
$$|\phi|\leq C|x|^{-\sigma}\|\phi\|_{*},\
\mbox{for\ all}\ |x|\leq 1,$$ and therefore, We obtain
\begin{equation}
\sup_{|x|\leq 1}|x|^{2+\sigma}w(x)^{p-2}|\phi|^{2} \leq
\|\phi\|_{*}^{2}\sup_{|x|\leq1}|x|^{2-\sigma}\\
=C\|\phi\|_{*}^{2}.
\end{equation}
From (3.3) and (3.4) it follows that \begin{equation}
\|w^{p-2}\phi^{2}\|_{**}\leq C\|\phi\|_{*}^{2}.\end{equation} To
estimate$\||\phi|^{p}\|_{**}$ we compute
\begin{equation}
\sup_{|x|\leq 1}|x|^{2+\sigma}|\phi(x)|^{p} \leq C\|\phi\|_{*}^{p}.
\end{equation}
Similarly
\begin{equation}
\sup_{|x|\geq 1}|x|^{2+\frac{2}{p-1}}|\phi|^{p} \leq
\|\phi\|^{p}_{*}.
\end{equation}
From (3.6) and (3.7) it follows that
\begin{equation}
\||\phi|^{p}\|_{**}\leq C\|\phi\|_{*}^{p}.
\end{equation}
 By (3.5) and
(3.8) we have
\begin{equation}\|N(\phi)\|_{**}\leq
C(\|\phi\|_{*}^{2}+\|\phi\|_{*}^{p}).
\end{equation}

Case $p<2$. In this case $|N(\phi)|\leq C|\phi|^{p}$ and hence, if
$0<\sigma\leq\frac{2}{p-1},$
\begin{equation}
\sup_{|x|\leq 1}|x|^{2+\sigma}|\phi(x)|^{p} \leq C\|\phi\|_{*}^{p}
\end{equation}
Similarly
\begin{equation}
\sup_{|x|\geq 1}|x|^{2+\frac{2}{p-1}}|\phi|^{p} \leq
\|\phi\|^{p}_{*}.
\end{equation}
 From (3.10)
and (3.11) it follows that for any $1<p<2$ and
$0<\sigma\leq\frac{2}{p-1},$
\begin{equation}
\|N(\phi)\|_{**}\leq C\|\phi\|_{*}^{p}.
\end{equation}
 From (3.9) and
(3.12) we have
\begin{equation}
\|N(\phi)\|_{**}\leq C(\|\phi\|_{*}^{2}+\|\phi\|_{*}^{p}).
\end{equation}

Now, we estimate $\|f_{\lambda}(x)\|_{**}$ as follows:
\begin{equation}
\sup_{|x|\leq 1}|x|^{2+\sigma}|f_{\lambda}(x)| =\sup_{\lambda
R_{1}<|x|\leq 1}|x|^{2+\sigma}|f_{\lambda}(x)|\leq
R_{1}^{\mu-2-\sigma}\to 0\ \ \mbox{as}\ \ R_{1}\to +\infty,
\end{equation}
provided that $\sigma=\frac{2}{p-1}$.

 Similarly,
\begin{equation}
\sup_{|x|\geq1}|x|^{2+\frac{2}{p-1}}|f_{\lambda}(x)|=o(\lambda).
\end{equation}

So we have, as $\lambda\to 0$,
\begin{equation}\|f_{\lambda}(x)\|_{**}=\sup_{|x|\leq
1}|x|^{2+\sigma}|f_{\lambda}(x)|+ \sup_{|x|\geq
1}|x|^{2+\frac{2}{p-1}}|f_{\lambda}(x)|\leq R_{1}^{\mu-2-\sigma}.
\end{equation}

We have already observed that $u=w+\phi$ is a solution of (1.5) if
$\phi$ satisfies equation (3.2). Consider the set
$$F=\{\phi \in L^{\infty}/\  \|\phi\|_{*}\leq \rho  \}$$
where $\rho \in (0,1)$ is to be chosen (suitably small) and the
operator$$\hbar(\phi)=T(N(\phi)-f_{\lambda}(x)).$$ We now prove that
$\hbar$ has a fixed point in $F$. For $\phi\in F$ we have
\begin{eqnarray*}
 \|\hbar(\phi)\|_{*}&\leq
&C\|N(\phi)\|_{**}+C\|f\|_{**}\\
&\leq&
C(\|\phi\|_{*}^{2}+\lambda^{-2}\|\phi\|_{*}^{p}+R_{1}^{\mu-2-\sigma})
\end{eqnarray*}
by (3.9) and (3.10), if  $\sigma=\frac{2}{p-1}$. Then we have
$$\|\hbar(\phi)\|_{*,\lambda}\leq
C(\rho^{2}+\rho^{p}+R_{1}^{\mu-2-\sigma})\leq \rho,$$ if we choose
$\rho$ is small enough and $R_{1}$ is large enough. Hence
$\hbar(F)\subset F$.

Now we show that $\hbar$ is a contraction mapping in $F$. Let us
take $\phi_{1},\phi_{2}$ in $F$, then
\begin{equation}
\|\hbar(\phi_{1})-\hbar(\phi_{2})\|_{*}\leq
C\|N(\phi_{1})-N(\phi_{2})\|_{**}.
\end{equation} Write
$$N(\phi_{1})-N(\phi_{2})=D_{\phi}N(\bar{\phi})(\phi_{1}-\phi_{2}),$$
where $\bar{\phi}$ lies in the segment joining $\phi_{1}$ and
$\phi_{2}$.

For $|x|\leq 1$,
$$ |x|^{2+\sigma}|N(\phi_{1})-N(\phi_{2})|
\leq |x|^{2}|D_{\phi}N(\bar{\phi})|\|\phi_{1}-\phi_{2}\|_{*},$$
while, for $|x|\geq 1$,
$$|x|^{2+\frac{2}{p-1}}|N(\phi_{1})-N(\phi_{2})| \leq
|x|^{2}|D_{\phi}N(\bar{\phi})|\|\phi_{1}-\phi_{2}\|_{*}.$$ Then we
have
\begin{equation}
\|N(\phi_{1})-N(\phi_{2})\|_{**}\leq
 C\sup_{x}|x|^{2}|D_{\phi}N(\bar{\phi})|\|\phi_{1}-\phi_{2}\|_{*}.
 \end{equation}
Directly from the definition of $N$, we compute
\begin{equation}
D_{\phi}N(\bar{\phi})=p[(w+\bar{\phi})^{p-1}-w^{p-1}].
\end{equation}
Thus $$|D_{\phi}N(\bar{\phi})|\leq
C(w^{p-2}|\bar{\phi}|+|\bar{\phi}|^{p-1}).$$ For all $x$ we have
\begin{equation}
|x|^{2}w^{p-2}|\bar{\phi}|\leq
C(\|\phi_{1}\|_{*}+\|\phi_{2}\|_{*})\leq C\rho.
\end{equation}
Similarly, for all $x$
\begin{equation}|x|^{2}|\bar{\phi}|^{p-1}\leq
C(\|\phi_{1}\|_{*}^{p-1}+\|\phi_{2}\|_{*}^{p-1})\leq
C\rho^{p-1}.
\end{equation}
 Estimates (3.13)-(3.15) show that
\begin{equation} \sup_{x}|x|^{2}|D_{\phi}N(\bar{\phi})|\leq
C(\rho+\rho^{p-1}).\end{equation} Gathering relations (3.11),
(3.12), (3.16) we conclude that $\hbar$ is a contraction mapping in
$F$, and hence a fixed point in this region indeed exists. So
$w+\phi_{\lambda}$ is solution of
\begin{equation}  \Delta u +u^{p}+f_{\lambda}(x)=0, \ \ \mbox{in}\ \ R^{n} \ \  u>0,
\ \ \lim_{|x|\to+\infty}u(x)=0, \end{equation} and
$$\phi_{\lambda}(x)\leq C\ \ \ \mbox{for all} \ \ \ x \in R^{n}\setminus \{0\}.$$
Thus $u_{\lambda}(x)=\lambda^{\frac{2}{p-1}}(w(\lambda
x)+\phi_{\lambda}(\lambda x))$ is a continuum of (1.1) and
$$\lim_{\lambda \to 0} u_{\lambda}(x)=0$$ uniformly in
$R^{n}\setminus \{0\}$. This finishes the proof of the theorem 1.\ \ $\Box$ \vskip 0.2in

{\it  To conclude, in this paper, instead of using sub-super solution method (which limits the
applicability on the exponent $p$), we use asymptotic analysis and
Liapunov-Schmidt reduction method to solve a open problem.}

\end{document}